\documentclass[12pt,reqno]{amsart}

\usepackage{amssymb,latexsym, amsmath, amscd, array}

\numberwithin{equation}{section}
\numberwithin{figure}{section}
\numberwithin{table}{section}

\DeclareMathOperator{\length}{{\rm length}}

\DeclareMathOperator{\area}{{\rm area}}

\DeclareMathOperator{\card}{{\rm card}}

\DeclareMathOperator{\sys}{{\rm sys} \pi_1}

\DeclareMathOperator{\SR}{{\rm SR}}

\def\ov{\overline}
\def\wh{\widehat}

\DeclareMathOperator{\Int}{{\rm Int}}
\DeclareMathOperator{\NN}{{\mathcal N}}

\def\rp{\mathbb R\mathbb P}

  \def\ov{\overline}

\newcommand{\R}{\mathbb R}

\newcommand{\Z}{\mathbb Z}
\newcommand{\gmetric}{{\mathcal G}}

\newtheorem{theorem}{Theorem}[section]

\newtheorem{proposition}[theorem]{Proposition}
\newtheorem{lemma}[theorem]{Lemma}

\newtheorem{prop}[theorem]{Proposition}

\newtheorem{cory}[theorem]{Corollary}
\newtheorem{claim}[theorem]{Claim}

\theoremstyle{definition}
\newtheorem{definition}[theorem]{Definition}
\newtheorem{example}[theorem]{Example}
\newtheorem{remark}[theorem]{Remark}

\newtheorem{question}[theorem]{Question}

\newcommand\theoref{Theorem~\ref}
\newcommand\lemref{Lemma~\ref}
\newcommand\propref{Proposition~\ref}
\newcommand\corref{Corollary~\ref}

\def\ie {{\it i.e.\ }}  
\def\cf {\hbox{\it cf.\ }}

\def\dist{{\rm dist}}

\def\ga{\alpha}

\def\eps{\varepsilon}

\def\m{\medskip}

\long\def\forget#1\forgotten{} %

\date{\today}
\begin{document}

\author[Y.~Rudyak]{Yuli B. Rudyak$^{*}$} 
\address{Department of Mathematics, University of Florida, PO Box 118105,
Gainesville, FL 32611-8105 USA} 
\email{rudyak@math.ufl.edu}
\thanks{$^{*}$Supported by NSF, grant 0406311}

\author[S.~Sabourau]{St\'ephane Sabourau}
\address{Laboratoire de Math\'ematiques et Physique Th\'eorique, Universit\'e 
de
Tours, Parc de Grandmont, 37400 Tours, France}
\address{Department of Mathematics, University of Pennsylvania, 209 South 33rd 
Street, Philadelphia, PA 19104-6395, USA}
\email{sabourau@lmpt.univ-tours.fr}

\title[Grushko meets systole]{Systolic invariants of groups and 2-complexes via 
Grushko decomposition}



\begin{abstract}

We prove a finiteness result for the systolic area of groups,
answering a question of M.~Gromov.  Namely, we show that there are
only finitely many possible unfree factors of fundamental groups
of~$2$-complexes whose systolic area is uniformly bounded. 
Furthermore, we prove a uniform systolic inequality for all 2-complexes 
with unfree fundamental group that improves the previously known bounds in this 
dimension.
\end{abstract}

\maketitle

\tableofcontents

\section{Introduction}
\label{sec:intro}

Throughout the article the word ``complex'' means ``finite simplicial 
complex'', unless something else is said explicitly.

\m 
Consider a piecewise smooth metric~$\gmetric$ on a complex~$X$.
The systole of~$\gmetric$, denoted~$\sys(\gmetric)$,
is defined as the least length of a noncontractible loop in~$X$. 
We define the systolic ratio of~$\gmetric$ as
\begin{equation}
\SR(\gmetric) = \frac{\sys(\gmetric)^{2}}{\area(\gmetric)},
\end{equation}
and the systolic ratio of~$X$ as
\begin{equation}
\SR(X) = \sup_{\gmetric} \SR(\gmetric),
\end{equation}
where the supremum is taken over the space of all the piecewise flat
metrics~$\gmetric$ on~$X$.
Note that taking the supremum over the space of all piecewise smooth metrics 
on~$X$ would yield the same value, \cf \cite{az}, \cite[\S3]{bz}.

We also define the {\em systolic ratio} of a finitely presentable
group~$G$ as
\begin{equation}
\label{103}
\SR(G) = \sup_{X}\; \SR(X),
\end{equation}
where the supremum is taken over all finite~$2$-complexes~$X$ with fundamental 
group isomorphic to~$G$.  It is also
convenient to introduce the {\em systolic area\/} $\sigma(G)$ of~$G$,
\cf~\cite[p.~337]{Gr2}, by setting
$$
\sigma(G) = \SR(G)^{-1}.
$$
Similarly, we define the systolic area of a~$2$-complex~$X$ and
of a piecewise flat metric~$\gmetric$ on~$X$ as
$\sigma(X) = \SR(X)^{-1}$ and~$\sigma(\gmetric) = \SR(\gmetric)^{-1}$, 
respectively.

\m
In this article, we study the systolic ratio of groups, or equivalently the 
systolic ratio of $2$-complexes.
Before stating our results, let us review what was previously known on the 
subject.

\m M.~Gromov \cite[6.7.A]{Gr1} (note a
misprint in the exponent) showed that every $2$-complex~$X$ with unfree 
fundamental group satisfies
the systolic inequality
\begin{equation}
\label{eq:72}
\SR(X) \leq 10^{4}.
\end{equation}
Contrary to the case of surfaces, where a (better) systolic inequality
can be derived by simple techniques, the proof of
inequality~\eqref{eq:72} depends on the advanced filling techniques
of~\cite{Gr1}.

Recently, in collaboration with M.~Katz, we improved the bound~\eqref{eq:72} 
using ``elementary'' techniques and characterized the $2$-complexes satisfying 
a systolic inequality, \cf \cite{krs}.
Specifically, we showed that every $2$-complex~$X$ with unfree fundamental 
group 
satisfies
\begin{equation}
\label{eq:12}
\SR(X) \leq 12.
\end{equation}
Furthermore, we proved that $2$-complexes with unfree fundamental groups 
are 
the only $2$-complexes satisfying a systolic inequality,
\ie for which the systolic ratio is bounded, \cf \cite{krs}. 

If one restricts oneself to surfaces, numerous systolic inequalities are 
available.
These inequalities fall into two categories.
The best estimates for surfaces of low Euler characteristic can be found in 
\cite{pu,Gr1,bav,Gr3,KS1,KS3,bcik}.
Near-optimal asymptotic bounds for the systolic ratio of surfaces of 
large genus have been established in \cite{Gr1,bal,KS2,sab} and \cite{BS,KSV}.

We refer to the expository texts \cite{Gr2,Gr3,CK} and the reference therein 
for an account on higher-dimensional systolic inequalities and other related 
curvature-free inequalities. 

\m
In order to state our main results, we need to recall Grushko decomposition in 
group theory.
By Grushko's theorem~\cite{St1, O}, every finitely generated group~$G$
has a decomposition as a free product of subgroups
\begin{equation}
\label{41b}
G=F_{p}*H_{1}*\dots*H_{q}
\end{equation}
such that~$F_{p}$ is free of rank~$p$, while every~$H_{i}$ is
nontrivial, non isomorphic to~$\Z$ and freely indecomposable.
Furthermore, given another decomposition of this sort, say
$G=F_{r}*H'_{1}*\dots*H'_{s}$, one necessarily has~$r=p$, $s=q$ and,
after reordering,~$H'_{i}$ is conjugate to~$H_{i}$.
We will refer to the number~$p$ in decomposition \eqref{41b} as the
{\em Grushko free index\/} of~$G$.

Thus, every finitely generated group~$G$ of Grushko free index~$p$ can
be decomposed as
\begin{equation}
G = F_{p} * H_{G},
\end{equation}
where~$F_{p}$ is free of rank~$p$ and~$H_{G}$ is of zero Grushko free
index.  The subgroup~$H_{G}$ is unique up to isomorphism.  Its
isomorphism class is called the {\em unfree factor} of (the isomorphism
class of)~$G$.

The Grushko free index of a complex is defined as the
Grushko free index of its fundamental group. 

\m
One of our main results answers, to a certain extent,
a question of M.~Gromov \cite[p.~337]{Gr2} about the systolic ratio of
groups.  More precisely, we obtain the following finiteness result.

\begin{theorem} \label{theo:fin}
Let~$C >0$.  The isomorphism classes of the unfree factors of the
finitely presentable groups~$G$ with $\sigma(G) < C$
lie in a finite set with at most
\[
A^{C^{3}}
\]
elements, where~$A$ is an explicit universal constant. 
\end{theorem}

While proving Theorem~\ref{theo:fin}, we improve the systolic inequalities 
\eqref{eq:72} and~\eqref{eq:12}.

\begin{theorem} \label{theo:main}
Every unfree finitely presentable group satisfies the inequality
\begin{equation} 
\label{eq:intro1}
\SR(G) \leq \frac{16}{\pi}.
\end{equation}
\end{theorem}

It is an open question if every $2$-complex with unfree fundamental group 
satisfies Pu's
inequality for~$\rp^{2}$, equivalently if the optimal constant
in~\eqref{eq:intro1} is~$\frac{\pi}{2}$. 

\m
The article is organized as follows.  
In Section~\ref{prelim}, we recall some topological preliminaries.
In Section~\ref{sec:pointed}, we investigate the geometry of pointed systoles 
and 
establish a lower bound on the area of ``small'' balls on $2$-complexes with 
zero Grushko free index. This yields a systolic inequality.
The existence of ``almost extremal regular'' metrics is established in
Section~\ref{sec:reg}.
Section~\ref{sec:counting} contains some combinatorial
results: we count the number of fundamental groups of complexes with some 
prescribed properties.
Using these results, we derive two finiteness results about the
fundamental groups of certain~$2$-complexes in Section~\ref{sec:finite}.
In Section~\ref{sec:dicho}, we relate the systolic ratio of a group to the 
systolic ratio of the free product of this group with~$\Z$.
In the last section, we combine all the results from the previous sections to 
prove our main theorems.

\m

{\bf Acknowledgment.}
The authors are very much indebted to Misha Katz for numerous exchanges during 
the preparation of this article.

\section{Topological preliminaries}
\label{prelim}

A proof of the following result, derived easily from the Seifert--van Kampen 
Theorem, can be found in~\cite{krs}.
 
\begin{lemma}\label{pi-exc}
Let~$(X,A)$ be a~$CW$-pair with~$X$ and~$A$ connected.  
If the inclusion~$j: A \to X$ induces the zero
homomorphism $j_*: \pi_1(A) \to \pi_1(X)$ of fundamental groups, then
the quotient map~$q: X \to X/A$ induces an isomorphism of fundamental
groups.
\end{lemma}

Let~$X$ be a finite connected complex and let~$f: X \to \R$ be a
function on~$X$.  Let
\[
[f\le r]:=\{x\in X\bigm| f(x)\le r\} \text{ and } [f\ge r]:=\{x\in X\bigm| 
f(x)\ge r\}
\]
denote the sublevel and superlevel sets of~$f$, respectively.

\begin{definition} \label{def:coalesce}
Suppose that a single path-connected component
of the superlevel set~$[f\ge r]$ contains
$k$ path-connected components of the level set~$f^{-1}(r)$.
Then we will say that the~$k$ path-connected components {\em
coalesce forward}.
\end{definition}

We will need the following result (we refer to~\cite{krs} for a more complete 
statement and a more detailed proof).

\begin{lemma} 
\label{42}
Assume that the pairs~$([f\ge r],f^{-1}(r))$ and~$(X,[f\leq r])$ are
homeomorphic to~$CW$-pair. Suppose that the set~$[f\le r]$ is
connected and that two connected components of~$f^{-1}(r)$ coalesce
forward. 
If the inclusion
\[
[f\le r]\subset X
\]
of the sublevel set~$[f \leq r]$ induces the zero homomorphism of
fundamental groups, then the Grushko free index of~$X$ is positive.
\end{lemma}

\begin{proof}Let $Y=[f \geq r]/\sim$ where~$x\sim y$ if
and only if~$x,y$ belong to the same component of~$f^{-1}(r)$.
The images~$a_i$ of the components of~$f^{-1}(r)$ under the quotient map $[f 
\geq r] \to Y$ form a finite set~$A \subset Y$.
By assumption, two points of~$A$ are joined by an arc in~$[f \geq r]$.
Therefore, the space~$Y \cup CA$, obtained by gluing an abstract cone over~$A$ 
to~$Y$, is homotopy equivalent to the wedge of~$S^1$ with another space~$Z$.
Hence,
\[
X/[f\le r]= Y/A \simeq Y \cup CA \simeq S^1 \vee Z.
\]
Thus, by the Seifert--van Kampen Theorem, the Grushko
free index of~$\pi_1(X/[f\le r])$ is positive.
Since the inclusion~$[f\le r]\subset X$ induces
the zero homomorphism of fundamental groups, we conclude that the group 
\mbox{$\pi_1(X/[f\le r])$}
is isomorphic to~$\pi_1(X)$ by \lemref{pi-exc}.
\end{proof}

We will also need the following technical result.

\begin{proposition} \label{homeo}
A level set of the distance function~$f$ from a point in a piecewise flat 
$2$-complex~$X$ is a finite graph.
In particular, the triangulation of~$X$ can be refined in such a way that the 
sets
\mbox{$[f \leq r]$}, $f^{-1}(r)$ and~$[f \geq r]$ become $CW$-subspaces
of~$X$.

Furthermore, the function~$\ell(r) = \length f^{-1}(r)$ is piecewise continuous.
\end{proposition}

Proposition~\ref{homeo} is a consequence of standard results in real algebraic 
geometry, \cf \cite{bcr}.
Indeed, note that $X$ can be embedded into some~$\R^{N}$ as a semialgebraic set 
and that the distance function~$f$ is a continuous semialgebraic function 
on~$X$.
Thus, the level curve~$f^{-1}(r)$ is a semialgebraic subset of~$X$ and, 
therefore, a finite graph, \cf \cite[\S9.2]{bcr}.
A more precise description of the level curves of~$f$ appears in~\cite{krs}.

The second part of the proposition also follows from~\cite[\S9.3]{bcr}.

\section{Complexes of zero Grushko free index} 
\label{sec:pointed}

The results of this section will be used repeatedly in the sequel.
These results also appear in~\cite{krs}.
We duplicate them here for the reader's convenience.

\begin{definition}
\label{df:ball}
Let~$X$ be a complex equipped with a piecewise smooth metric.  A
shortest noncontractible loop of~$X$ based at~$x \in X$ is called a
{\it pointed systolic loop at~$x$. \/} Its length, denoted
by~$\sys(X,x)$, is called the {\it pointed systole at~$x$.\/} 
\end{definition}

\m As usual, given $x\in X$ and $r\in \R$, we denote by $B(x,r)$ the ball of 
radius $r$ centered at $x$, $B(x,r)=\{a\in X\bigm|\dist(x,a) \leq r\}$.
 
\begin{prop}\label{ball}
If $r < \frac{1}{2} \sys(X,x)$ then the inclusion~$B(x, r)\subset X$ induces
the zero homomorphism of fundamental groups.
\end{prop}

\begin{proof}
Suppose the contrary.
Consider all the loops of $B(x,r)$ based at~$x$ that are noncontractible in~$X$. 
Let $\gamma \subset B(x,r)$ be the shortest of these loops.
We have $L=\length(\gamma) \geq \sys(X,x)$.
Let~$a$
be the point of~$\gamma$ that divides~$\gamma$ into two arcs
$\gamma_1$ and~$\gamma_2$ of the same length~$L/2$. Consider a shortest geodesic 
path~$c$, of length $d=d(x,a) < r$, that joins~$x$ to~$a$.
Since at least one of the curves~$\gamma_1\cup c_-$ or~$c\cup
\gamma_2$ is noncontractible, we conclude that~$d+L/2\ge L$, \ie~$d\ge
L/2$ (here~$c_-$ denotes the path~$c$ with the opposite
orientation). Thus
$$
\sys(X,x)>2r\ge 2d\ge L \ge\sys(X,x).
$$
That is a contradiction.
\end{proof} 

\forget
\begin{remark}
\label{ball2}
Alternatively,~$\sys(X,x)$ could be defined as twice the upper bound
of the reals~$R > 0$ such that the homomorphism $\pi_1(B(x,R))\to \pi_1(X)$ 
induced by the inclusion $B(x,R) \subset X$ is trivial.
In other words, every loop (not necessarily passing through~$x$) in 
$B(x,R)$ is contractible in~$X$.
\end{remark}
\forgotten

The following lemma describes the structure of a pointed systolic
loop.

\begin{lemma} 
\label{lem:simple}
Let $X$ be a complex equipped with a piecewise flat metric.
Let $\gamma$ be a pointed systolic loop at~$x \in X$ of length
$L=\sys(X,x)$.
\begin{enumerate}
\item[(i)] The loop~$\gamma$ is formed of two distance-minimizing
arcs, starting at~$p$ and ending at a common endpoint, of length
$L/2$.
\item[(ii)] Any point of self-intersection of the loop~$\gamma$ is no
further than $\tfrac{1}{2} \left(\sys(X,x) - \sys(X) \right)$ from~$x$.
\end{enumerate}
\end{lemma}

\begin{proof}
Consider the arc length parameterization $\gamma(s)$ of the loop $\gamma$ with
\mbox{$\gamma(0)=\gamma(L)=x$}.  Let~$y = \gamma \left(L/2 \right)
\in X$ be the ``midpoint'' of $\gamma$.  Then $y$ splits~$\gamma$ into
a pair of paths of the same length~$L/2$, joining~$x$ to~$y$.
By \propref{ball}, if $y$ were contained in the open ball
$B(x,L/2)$, the loop $\gamma$ would be contractible.  This
proves item~(i).

If $x'$ is a self-intersection point of~$\gamma$, the loop~$\gamma$
decomposes into two loops $\gamma_{1}$ and $\gamma_{2}$ based at~$x'$,
with $x \in \gamma_{1}$.  Since the loop~$\gamma_{1}$ is shorter than
the pointed systolic loop~$\gamma$ at~$x$, it must be contractible.
Hence~$\gamma_{2}$ is noncontractible, so that
\[
\length(\gamma_{2}) \geq \sys(X).
\]
Therefore,
\[
\length(\gamma_{1}) = L - \length(\gamma_2) \leq \sys(X,x) - \sys(X),
\]
proving item (ii).
\end{proof}

The following proposition provides a lower bound for the length of
level curves in a $2$-complex.

\begin{proposition} 
\label{lem:1compo}
Let $X$ be a piecewise flat $2$-complex.
Fix~$x \in X$.
Let $r$ be a real number satisfying
\[
\sys(X,x) - \sys(X) <2r< \sys(X,x) .
\]
Consider the curve~$S=\{a\in X\bigm|\dist(x,a)=r\}$.  Let~$\gamma$ be a pointed 
systolic loop at~$x$. 
If $\gamma$ intersects exactly one connected component of~$S$, then
\begin{equation} \label{eq:lengthS}
\length S \geq 2r - \left( \sys(X,x) - \sys(X) \right).
\end{equation}
\end{proposition}

\begin{proof}
By Lemma~\ref{lem:simple}, the loop~$\gamma$ is formed of two
distance-minimizing arcs which do not meet at distance~$r$ from~$x$.
Thus, the loop $\gamma$ intersects~$S$ at exactly two points.
Let~$\gamma' = \gamma \cap B$ be the subarc of $\gamma$ lying in~$B$.

If~$\gamma$ meets exactly one connected component of~$S$, there exists
an embedded arc $\alpha \subset S$ connecting the endpoints of
$\gamma'$.  By \propref{ball}, every loop based at~$x$ and lying
in~$B(x,r)$ is contractible in~$X$.  Hence~$\gamma'$ and $\alpha$ are
homotopic, and the loop $ \alpha \cup (\gamma \setminus \gamma')$ is
homotopic to~$\gamma$.  Hence,
\begin{equation}
\length(\alpha) + \length(\gamma) - \length(\gamma') \geq \sys(X).
\end{equation}
Meanwhile, $\length(\gamma) = \sys(X,x)$ and $\length(\gamma')=2r$, proving the
lower bound \eqref{eq:lengthS}, since~$\length(S) \geq
\length(\alpha)$.
\end{proof}

The following result provides a lower bound on the area of ``small'' balls of 
$2$-complexes with zero Grushko free index, \cf Section~\ref{sec:intro}.

\begin{theorem} \label{theo:smallballs}
Let~$X$ be a piecewise flat $2$-complex with zero Grushko free index.
Fix $x \in X$.
For every real number $R$ such that
\begin{equation}
\label{51z}
\sys(X,x)-\sys(X) \leq 2R \leq \sys(X,x),
\end{equation}
the area of the ball~$B(x,R)$ of radius~$R$ centered at~$x$ satisfies
\begin{equation} 
\label{eq:smallballs}
\area B(x,R) \geq \left( R- \tfrac{1}{2}(\sys(X,x)-\sys(X))
\right)^{2}.
\end{equation}
In particular, we have
\[
\SR(X) \leq 4.
\]
\end{theorem}

\begin{remark}
The example of a piecewise flat $2$-complex with a circle of length the 
systole of~$X$ attached to it shows that the assumption on the fundamental 
group 
of the complex cannot be dropped.
\end{remark}

\begin{proof}[Proof of Theorem~$\ref{theo:smallballs}$] 

Let $L= \sys(X,x)$.
Let $r$ be a real number satisfying $L - \sys(X) \leq 2r \leq L$.
Denote by~$S=S(x,r)$ and~$B=B(x,r)$, respectively,
the level curve and the ball of radius~$r$ centered at~$x$.
Let $\gamma$ be a pointed systolic loop at~$x$.

If $\gamma$ intersects two connected components of~$S$, then by 
Lemma~\ref{lem:simple}, there exists an arc of~$\gamma$ lying in~$X \setminus 
\Int(B)$, which joins these two components of~$S$.
That is, the components coalesce forward.
Thus, by Lemma~\ref{42} and \propref{ball}, the complex~$X$ has a positive 
Grushko free index, which is excluded.

Therefore, the loop~$\gamma$ meets a single connected component
of~$S$. Now, Proposition~\ref{lem:1compo} implies
that
\begin{equation}
\length S(r) \geq 2r - \left( \sys(X,x) - \sys(X) \right).
\end{equation}

Let $\eps = \sys(X,x) - \sys(X)$.
Using the coarea formula, \cf \cite[3.2.11]{Fe1}, \cite[13.4]{bz}, as in 
\cite[Theorem~5.3.1]{bz}, \cite{heb} and \cite[5.1.B]{Gr1}, we obtain
\begin{eqnarray*}
\area B(x,R) & \geq & \int_{\frac{\eps}{2}}^{R} \length S(r) \, dr \\
& \geq & \int_{\frac{\eps}{2}}^{R} (2r - \varepsilon) \, dr \\ & \geq 
&
\left( R-\frac{\eps}{2}\right)^{2}
\end{eqnarray*}
for every~$R$ satisfying \eqref{51z}.

Now, if we choose $x \in X$ such that a systolic loop passes through~$x$, then 
$\sys(X,x)=\sys(X)$.
In this case, setting $R = \tfrac{1}{2} \sys(X,x)$, we obtain $ \area(X) \geq 
\tfrac{1}{4} \sys(X)^{2}$, as
required.
\end{proof}

\section{Existence of~$\varepsilon$-regular metrics} \label{sec:reg}

\begin{definition}
A metric on a complex~$X$ is said to be
$\varepsilon$-regular if $\sys(X,x) < (1+\varepsilon) \sys(X)$ for
every~$x$ in~$X$.
\end{definition}

\begin{lemma} \label{lem:regular}
Let~$X$ be a $2$-complex with unfree fundamental group.  Given a
metric~$\gmetric$ on~$X$ and~$\varepsilon > 0$, there exists an
$\varepsilon$-regular piecewise flat metric~$\gmetric_{\varepsilon}$ on~$X$ 
with 
a
systolic ratio as good as for~$\gmetric$, \ie
$\SR(\gmetric_{\varepsilon}) \geq \SR(\gmetric)$.
\end{lemma}

\begin{proof}
We argue as in~\cite[$5.6.C''$]{Gr1}.  Choose~$\varepsilon' >0$ such
that~$\varepsilon' < \min\{\varepsilon,1\}$.
Fix~\mbox{$r'=\frac{1}{2} \varepsilon' \sys(\gmetric)$} and $r>0$,
with~$r<r'$.  Subdividing~$X$ if necessary, we can assume that the
diameter of the simplices of~$X$ is less than~$r'-r$.  The
{\em approximating ball}~$B'(x,r)$ is defined
as the union of all simplices of~$X$ intersecting~$B(x,r)$.  By 
construction,~$B'(x,r)$ is a
path connected subcomplex of~$X$ which contains~$B(x,r)$ and is
contained in $B(x,r')$.  In, particular,
the inclusion~$B' \subset X$ induces the trivial homomorphism of
fundamental groups.

Assume that the metric~$\gmetric_0=\gmetric$ on~$X_{0}=X$ is
not
already~$\varepsilon'$-regular.
There exists a point~$x_{0}$ of~$X_{0}$ such that
\begin{equation} \label{eq:nonreg}
\sys(X_{0},x_{0}) > (1+\varepsilon') \sys(X_{0}).
\end{equation}
Consider the space
$$
X_{1} = X_{0}/B'_{0}
$$
obtained by collapsing the approximating ball~$B'_{0}:=B'(x_0,r)$.
Denote by~$\gmetric_1$ the length structure induced
by~$\gmetric_0$ on~$X_1$.
Let~$p_{0}:X_{0} \longrightarrow X_{1}$ be the (non-expanding)
canonical projection.
By \lemref{pi-exc}, the projection~$p_{0}$ induces an
isomorphism of fundamental groups.
Consider a systolic loop~$\gamma$ of~$\gmetric_1$.
Clearly,~$\length_{\gmetric_1}(\gamma) \leq
\sys(\gmetric_0)$.

If~$\gamma$ does not pass through the point~$p_0(B'_0)$, then the
preimage of~$\gamma$ under~$p_0$ is a noncontractible loop of the same
length as~$\gamma$.
Therefore,~$\sys(\gmetric_1) = \sys(\gmetric_0)$.

Otherwise,~$\gamma$ is a loop based at the
point~$p_{0}(B'_{0})$.
It is possible to construct a (noncontractible) loop
$\overline{\gamma}$ on~$X_0$ passing through~$x_0$ with
$$
\length_{\gmetric_0}(\overline{\gamma}) \leq
\length_{\gmetric_1}(\gamma) + 2r',
$$
whose image under~$p_0$ agrees with~$\gamma$.
From~\eqref{eq:nonreg}, we deduce that
$$
\length_{\gmetric_1}(\gamma) \geq
\length_{\gmetric_0}(\overline{\gamma}) -2r' \geq
(1+\varepsilon')
\sys(\gmetric_0) - 2r' = \sys(\gmetric_0).
$$
Thus, the systole of~$\gmetric_1$ is the same as the systole
of~$\gmetric_0$ and its area (or Hausdorff measure)
is at most the area of~$\gmetric_0$. Hence,~$\SR(\gmetric_1)\ge\SR(\gmetric_0)$

If~$\gmetric_1$ is not~$\varepsilon'$-regular, we apply the
same
process to~$\gmetric_1$.
By induction, we construct:
\begin{itemize}
\item a sequence of points~$x_i \in X_i$ with
$$
\sys(X_{i},x_{i}) > (1+\varepsilon') \sys(X_{i}), 
$$
\item a sequence of approximating
balls~$B'_i:=B'(x_i,r)$ in~$X_i$, 

\item a sequence of spaces~$X_{i+1}$ obtained
from~$X_i$ by
collapsing~$B'_{i}$ into a point (with~$\pi_{1}(X_{i+1})
\simeq
\pi_{1}(X_{i})$), 

\item a sequence of metrics~$\gmetric_{i+1}$ induced
by~$\gmetric_i$ on~$X_{i+1}$, 

\item a sequence of canonical projections~$p_i:X_i
\longrightarrow
X_{i+1}$.

\end{itemize}
This process stops when we obtain an~$\varepsilon'$-regular
metric
(with a systolic ratio as good as the one of~$\gmetric$).

Now we show that this process really stops.
Let~$B_1^i,\dots,B_{N_i}^i$ be a maximal system of disjoint
balls of
radius~$r/3$ in~$X_i$.
Since~$p_{i-1}$ is non-expanding, the
preimage~$p_{i-1}^{-1}(B_k^i)$
of~$B_k^i$ contains a ball of radius~$r/3$ in~$X_{i-1}$.
Furthermore, the preimage~$p_{i-1}^{-1}(x_i)$ of~$x_i$
contains a
ball~$B_{i-1}$ of radius~$r$ in~$X_{i-1}$.
Thus, two balls of radius~$r/3$ lie in the preimage of~$x_i$
under~$p_{i-1}$.
It is then possible to construct a system of~$N_i+1$ disjoint
disks
of radius~$r/3$ in~$X_{i-1}$.
Thus,~$N_{i-1} \geq N_i +1$ where~$N_i$ is the maximal number
of
disjoint balls of radius~$r/3$ in~$X_i$.
Therefore, the process stops after~$N$ steps with~$N \leq
N_0$.

Let~$\pi=p_{N-1} \circ \dots \circ p_0$ be the projection
from~$X$
to~$X_N$.
Denote by~$\Delta$ the set formed of the points of~$X_{N}$
whose
preimage under~$\pi$ is a singleton, \ie
$$
\Delta = \{ y \in X_{N} \mid \card \pi^{-1}(y) = 1 \}.
$$
By construction, the set~$X_N \setminus \Delta$ has at
most~$N$
points, which will be called the singularities of~$X_N$.

Let~$\gmetric_{t}$ be the length structure on~$X$ induced by
$e^{-t\varphi} \gmetric$, where~$t>0$ and
$\varphi(x)=\dist_{\gmetric}(\pi^{-1}(\Delta),x)$ for $x\in X$.
of~$X$ 
Clearly,~$\area(\gmetric_{t}) \leq \area(\gmetric)$ and
$\sys(\gmetric_{t}) \geq \sys(X_{N}) = \sys(\gmetric)$.
Therefore,~$\SR(\gmetric_{t}) \geq \SR(\gmetric)$.

It suffices to prove that~$\gmetric_{t}$ is~$\varepsilon$-regular for~$t$
large enough. Since~$X_{N}$ is~$\varepsilon'$-regular, this follows, in turn,
from the claim~\ref{uniform} below.

Strictly speaking, the metrics~$\gmetric_{t}$ are not piecewise flat but we can 
approximate them by piecewise flat metrics as in~\cite{az} (see 
also~\cite[\S3]{bz}) to obtain the desired conclusion.

\begin{claim}\label{uniform}
The family~$\{\sys(\gmetric_{t},x)\}$ converges to~$\sys(X_{N},\pi(x))$ 
uniformly in~$x$ as $t$ goes to infinity.
\end{claim}

Clearly, for every~$x$ in~$X$ and~$t > 0$, we have
\begin{equation} \label{eq:ucv1}
\sys(X_{N},\pi(x)) \leq \sys(\gmetric_{t},x).
\end{equation}
Fix~$\delta >0$. Take a pointed systolic loop~$\gamma \subset X_N$ at some 
fixed point~$y$ of~$X_{N}$ and let $\gamma$ pass through~$k(y)$ singularities
of~$X_N$.
Given~$z \in X_N$ at distance at most~$R=\delta/5$ from~$y$,
the loop
$[z,y] \cup \gamma \cup [y,z]$ based at~$z$, where~$[a,b]$
represents
a segment joining~$a$ to~$b$, is freely homotopic to~$\gamma$
and
passes through at most~$k(y) + 2N$ singularities.
Moreover, its length is at most~$\sys(X_{N},y) + 2R \leq
\sys(X_{N},z) + 4R$ since~$\sys(X_{N},.)$ is~$2$-Lipschitz.

Let~$k = \max_i k(y_i) + 2N$ where the~$y_i$'s are the centers
of a
maximal system of disjoint balls of radius~$R/2$ in~$X_N$.
It is possible to construct for every~$z$ in~$X_N$ a
noncontractible
loop~$\gamma_z$ based at~$z$ passing through at most~$k$
singularities of length at most~$\sys(X_{N},z) + 4R$.

The preimages~$U_{i}$ under~$\pi: X \to X_N$ of the singularities
of~$X_{N}$ are path-connected.
Choose~$t$ large enough so that every pair of points
in~$U_{i}$ can
be joined by an arc of~$U_i$ of~$\gmetric_{t}$-length less
than some
fixed~$\eta > 0$ with~$\eta < R/k$.
Fix~$x \in X$.
Consider the loop~$\gamma=\gamma_z$ of~$X_N$ based
at~$z=\pi(x)$
previously defined.
There exists a noncontractible loop~$\overline{\gamma}
\subset X$
based at~$x$ of length
$$
\length_{\gmetric_{t}}(\overline{\gamma}) \leq
\length_{X_{N}}(\gamma) + k \eta
$$
whose image under~$\pi$ agrees with~$\gamma$.
Therefore,
$$
\sys(\gmetric_{t},x) \leq \sys(X_{N},\pi(x)) + 4R +R.
$$
Hence,
\begin{equation} \label{eq:ucv2}
\sys(\gmetric_{t},x) \leq \sys(X_{N},\pi(x)) + \delta.
\end{equation}
Since~$\sys(\gmetric_{t'},x) \leq \sys(\gmetric_{t},x)$ for
every~$t'
\geq t$, the inequalities \eqref{eq:ucv1} and \eqref{eq:ucv2}
lead to
the desired claim.

This concludes the proof of the Lemma~\ref{lem:regular}.
\end{proof}

\section{Counting fundamental groups} 
\label{sec:counting}

Let~$X$ be a complex endowed with a piecewise flat metric, Consider a finite 
covering~$\{B_i\}$ of~$X$ by open balls of
radius~$R=\frac{1}{6} \sys(X)$.  Denote by~$\NN$ the nerve of this covering.

\begin{lemma} \label{lem:nerve}
The fundamental groups of~$X$ and~$\NN$ are isomorphic.
\end{lemma}

\begin{proof}
Recall that, by definition, the vertices~$p_i$ of~$\NN$ are identified
with the balls~$B_i$. Furthermore,~$k+1$ vertices~$p_{i_0}, \ldots,
p_{i_k}$ form a $k$-simplex of $\NN$ if and only if $B_{i_0}\cap \cdots \cap
B_{i_k}\ne \emptyset$.  Given~$x$ and~$y$ in~$X$, fix a minimizing
path (not necessarily unique), denoted by~$[x,y]$, from~$x$ to~$y$.

We denote by~$\NN_i$ the~$i$-skeleton of~$\NN$. Define a map
$f: \NN_{1} \longrightarrow X$ as follows.
The map~$f$ sends the vertices~$p_i$ to the centers~$x_i$ of
the balls~$B_i$ and the edges~$[p_i,p_j]$ to the
segments~$[x_i,x_j]$ (previously chosen).
By construction, the distance between two centers~$x_i$
and~$x_j$ corresponding to a pair of adjacent vertices is less
than~$2R$.
Thus, the map~$f$ sends the boundary of each~$2$-simplex
of~$\NN$ to
loops of length less than~$6R = \sys(X)$.
By definition of the systole, these loops are contractible
in~$X$.
Therefore, the map~$f$ extends to
a map~$F: \NN_{2} \longrightarrow X$.

Choose a center $x_{\ga(0)}$ of some of the balls $B_i$. We claim that the 
homomorphism~$F_*: \pi_1(\NN_{2}, p_{\ga(0)}) \longrightarrow
\pi_{1}(X, x_{\ga(0)})$ is an isomorphism. Since the nerve~$\NN$ and 
its~$2$-skeleton~$\NN_{2}$ have the
same
fundamental group, we conclude that~$\pi_{1}(X)$ and
$\pi_{1}(\NN)$ are isomorphic.

\m We prove the surjectivity of~$F_*$ only. The
injectivity can be proved in a similar way, we leave it to the reader.

\m
Given a piecewise smooth path~$\gamma:I \longrightarrow X$,  
$\gamma(0)=\gamma(1)=x_{\ga(0)}$, we construct the following path
$\overline{\gamma}:I \longrightarrow \NN_{1}, 
\ov\gamma(0)=\ov\gamma(1)=p_{\ga(0)}$ such that the loop~$F(\ov\gamma)$ is
homotopic to~$\gamma$. Fix a subdivision~$t_0=0<t_1 < \cdots < t_m<t_{m+1}=1$ 
of~$I$ such
that~$\gamma([t_{k},t_{k+1}])$ is contained in
some~$B_{\alpha(k)}$ and
the length of~$\gamma_{|[t_k,t_{k+1}]}$ is less than~$\frac{1}{3}$ for $k=0, 
\ldots, m$. The map~$\overline{\gamma}$ takes the segment~$[t_k,t_{k+1}]$
to the
edge~$[p_{\alpha(k)},p_{\alpha(k+1)}]$ of~$\NN$.
By construction, we have~$\overline{\gamma}(t_k) = p_{\alpha(k)}$ and
$F(\overline{\gamma}(t_k)) = x_{\alpha(k)}$.
Therefore, the image of~$\overline{\gamma}$ under $F$ is a
piecewise linear loop which agrees with the union
$$
\bigcup_{k=0}^m [x_{\alpha(k)},x_{\alpha(k+1)}]
$$
where the segments~$[x_{\alpha(k)},x_{\alpha(k+1)}]$ are
previously
fixed.

Consider the following loops of~$X$
$$
c_k = \gamma([t_k,t_{k+1}]) \cup
[\gamma(t_{k+1}),x_{\alpha(k+1)}]
\cup
[x_{\alpha(k+1)},x_{\alpha(k)}] \cup
[x_{\alpha(k)},\gamma(t_k)]
$$
where~$[x_{\alpha(k+1)},x_{\alpha(k)}]$ agrees with~$F \circ
\overline{\gamma}([t_k,t_{k+1}])$.
The length of~$c_k$ is
$$
\length(c_k) < \frac{1}{3} \sys(X) + R + 2R + R = \sys(X).
$$
Hence, the loop~$c_k$ is contractible.
Therefore, the loops~$\gamma$ and~$F\circ\overline{\gamma}$ are
homotopic, and thus the homomorphism~$F_*$ is
surjective.
\end{proof}

\begin{definition} \label{def:Gamma}
The isomorphism classes of the fundamental groups of the finite
$2$-complexes with at most~$n$ vertices form a
finite set~$\Gamma(n)$.  We define~$\Gamma'(n)$ as the union
of~$\Gamma(n)$ and the set formed of the unfree factors of the
elements of~$\Gamma(n)$.
\end{definition}

\begin{cory}\label{cor:counting}
Suppose that the covering $\{B_i\}$ of $X$ in Lemma $\ref{lem:nerve}$ consists 
of $m$ elements. Then $\pi_1(X)\in \Gamma(m)$.
\end{cory}

\begin{proof} This follows from \lemref{lem:nerve}, since the nerve of the 
covering has $m$ vertices.
\end{proof}

Now we estimate the numbers $\Gamma(n)$ and $\Gamma'(n)$.

\begin{lemma} \label{lem:count}
Up to isomorphism, the number of~$2$-dimensional simplicial
complexes
having~$n$ vertices is at most
\begin{equation}
2^{\frac{(n-1)n(n+1)}{6}} < 2^{n^3}.
\end{equation}
In particular, the sets~$\Gamma(n)$ and~$\Gamma'(n)$ contain
less
than~$2^{n^{3}}$ elements.
\end{lemma}

\begin{proof}
The maximal number of edges in a simplicial complex with~$n$
vertices
is
equal to the cardinality of~$\{(i,j) \mid 1 \leq i < j \leq n
\}$, which
is~$\frac{n(n-1)}{2}$.
Similarly, the maximal number of triangles in a simplicial
complex
with~$n$ vertices is equal to the cardinality of~$\{(i,j,k) \mid
1 \leq
i < j
<k \leq n \}$, which is~$\frac{n(n-1)(n-2)}{6}$.
Thus, the number of isomorphism classes of~$1$-dimensional
simplicial
complexes having~$n$ vertices is at most
\begin{equation} \label{eq:*}
2^{\frac{n(n-1)}{2}}.
\end{equation}
Therefore, the number of~$2$-dimensional simplicial complexes
whose
$1$-skeleton agrees with one of these~$1$-dimensional
simplicial
complexes
is at most
\begin{equation} \label{eq:**}
2^{\frac{n(n-1)(n-2)}{6}}.
\end{equation}
The product of~\eqref{eq:*} and~\eqref{eq:**} yields an upper
bound
on the
number of isomorphism classes of~$2$-dimensional simplicial
complexes
having~$n$ vertices.

Note that~$\Gamma'(n)$ has at most twice as many elements
as~$\Gamma(n)$.
The second part of the lemma follows then from the first part.
\end{proof}

\section{Two systolic finiteness results}
\label{sec:finite}

\begin{proposition} \label{prop:finite}
Let~$X$ be a $2$-complex equipped with a piecewise flat metric.
Suppose that the area of every ball~$B(R)$ of
radius~$R=\frac{1}{12} \sys(X)$ in~$X$ is at least~$\alpha
\sys(X)^{2}$, \ie
\begin{equation}\label{alpha}
\area B(R) \geq \alpha \sys(X)^{2}.
\end{equation}
If~$\sigma(X) < C$, then the isomorphism class of
the
fundamental group of~$X$ lies in the finite
set~$\Gamma(C/\alpha)$.
\end{proposition}

\begin{proof}
Consider a maximal system of disjoint open balls~$B(x_{i},R)$
in~$X$ of radius~$R=\frac{1}{12}\sys(X)$ with  centers~$x_{i}, i=1, \ldots, m$. 
By the assumption,
\begin{equation}
\area B(x_{i},R) \geq \alpha \sys(X)^{2}.
\end{equation}
Therefore, this system admits at most~$\frac{\area(X)}{\alpha
\sys(X)^{2}}$ balls.
Thus,
\begin{equation} \label{eq:I}
m \leq C/\alpha.
\end{equation}
The open balls~$B_{i}$ of radius~$2R=\frac{1}{6}\sys(X)$
centered
at~$x_{i}$ form a covering of~$X$.
From \corref{cor:counting}, the fundamental group of~$X$ lies
in~$\Gamma(m)\subset \Gamma(C/\alpha)$.
\end{proof}

\begin{theorem} \label{theo:finite}
Given~$C > 0$, there are finitely many isomorphism classes of finitely
presented groups~$G$ of zero Grushko free index  such that \mbox{$\sigma(G) < 
C$}.

More precisely, the isomorphism class of every finitely
presented group~$G$ with zero Grushko free index and~$\sigma(G) < C$ lies in 
the finite set~$\Gamma(144 \, C)$, which has at most
$$ K^{C^{3}},$$
elements. Here, $K$ is an explicit universal constant.
\end{theorem}

\begin{remark}
Clearly, we have~$\sigma(G_{1}*G_{2}) \leq \sigma(G_{1}) +
\sigma(G_{2})$ for every finitely presentable groups~$G_{1}$ and
$G_{2}$ (by taking the wedge of corresponding complexes).
In particular, the inequality~$\sigma(F_{p}*G) \leq \sigma(G)$ holds for
every~$p$.  Thus, the assumption that~$G$ has zero Grushko free index  cannot 
be
dropped in the previous finiteness result.
\end{remark}

\begin{question}

For which groups, $G_{1}$ and~$G_{2}$, does the relation $\sigma(G_{1}*G_{2}) = 
\sigma(G_{1}) + \sigma(G_{2})$ hold?
\end{question}

\begin{proof}[Proof of Theorem~$\ref{theo:finite}$]
Consider a finitely presentable group $G$ of zero Grushko free index  and
such that~$\sigma(G) < C$.
There exist a $2$-complex~$X$ with
fundamental
group
isomorphic to~$G$ and a piecewise flat metric~$\gmetric$ on~$X$ such that 
$\sigma(\gmetric) < 
C$.
Let~$0 < \varepsilon < \frac{1}{12}$.
Fix a~$2 \varepsilon$-regular piecewise flat metric on~$X$ with a better
systolic
ratio
than the one of~$\gmetric$, \cf Lemma~\ref{lem:regular}.
By \theoref{theo:smallballs},
\begin{equation}
\area B(R) \geq \left(\frac{1}{12} - \varepsilon\right)^{2}
\sys(X)^{2}
\end{equation}
for all balls~$B(R)$ of radius~$R=\frac{1}{12}\sys(X)$.
Since~$\sigma(X) < C$, we deduce from
Proposition~\ref{prop:finite} that the isomorphism class of
the fundamental group of~$X$ lies in the finite set
$$
\Gamma
\left(
\frac{C}{(\frac{1}{12}-\varepsilon)^{2}}
\right)=\Gamma\left(\frac{144\,C}{(1-12\varepsilon)^2}\right)
$$
for every $\varepsilon >0$ small enough.
Thus, the isomorphism class of~$G$ lies in~$\Gamma(144 \, C)$.

By Lemma~\ref{lem:count}, this set has at most
$$
(2^{12^6})^{C^{3}}
$$
elements. Hence the result.
\end{proof}

\begin{example}
\label{moore}
It follows from \theoref{theo:finite} that the systolic ratio of the
cyclic groups~$\Z/n\Z$ of order~$n$ goes to zero as~$n \to \infty$,
\ie
$$
\lim_{n \rightarrow \infty} \SR(\Z/n\Z) = 0.
$$ It would be interesting, however, to evaluate the
value~$\SR(\Z/n\Z)$. 
\end{example}

\section{Systolic area comparison}
\label{sec:dicho}

\m
Let $G$ be an unfree finitely presentable group with $G = F_{p}*H$
where $F_{p}$ is free of rank~$p$ and $H$ is of zero Grushko free
index.
Fix~$\delta \in (0,\frac{1}{12})$ (close to zero) and~$\lambda
> \frac{1}{\pi}$ (close to~$\frac{1}{\pi}$).  Choose $\varepsilon <
\delta$ (close to zero) such that~$0< \varepsilon < 4 (\lambda
-\frac{1}{\pi}) (\delta-\varepsilon)^{2}$.  From
Lemma~\ref{lem:regular}, there exists a $2$-complex~$X$ with fundamental 
group isomorphic to~$G$ and a $2 \varepsilon$-regular
piecewise flat metric~$\mathcal{G}$ on~$X$ such that
\begin{equation}\label{inequality}
\sigma(\mathcal{G}) \leq \sigma(G) +\varepsilon.
\end{equation}
We normalize the metric~$\gmetric$ on~$X$ so that its systole
is equal to~$1$.

Denote by~$B(x,r)$ and~$S(x,r)$ the ball and the sphere of
radius~$r <
\frac{1}{2}$ centered at some point~$x$ of~$X$.
Note that 
\begin{equation} \label{eq:delta}
\delta > \varepsilon > \frac{1}{2} \left( \sys(X,x) - \sys(X) \right)
\end{equation}
for every $x \in X$.

\begin{lemma}\label{path} Suppose that there exist
$x_0 \in X$ and~$r_0 \in (\delta, \frac{1}{2})$ such that
\begin{equation} \label{eq:assum}
\area B > \lambda \, (\length S)^2
\end{equation}
where~$B=B(x_0,r_0)$ and~$S=S(x_0,r_0)$.
Then, the Grushko free index ~$p$ of~$G$ is positive, and
\begin{equation}
\sigma(G) \geq \sigma(F_{p-1}*H) - \varepsilon.
\end{equation}
\end{lemma}

\begin{proof}
First, we prove that $p>0$. We let~$f(x)=\dist(x_{0},x)$ and show that two 
path-connected components of~$S=f^{-1}(r_{0})$
coalesce forward, \cf 
Definition~\ref{def:coalesce} and \lemref{42}.  Denote by~$\overline{X}$ the
$2$-complex obtained from~$X$ by attaching
cones~$C_{i}$ over each connected component~$S_{i}$ of~$S$,~$1 \leq i
\leq m$.  By \propref{ball}, the connected components~$S_{i}$
are contractible in~$X$.  Therefore, the fundamental groups of~$X$
and~$\overline{X}$ are isomorphic, \ie
\begin{equation} \label{eq:isom}
\pi_{1}(X) \simeq \pi_{1}(\overline{X}).
\end{equation}
Fix a segment~$[x_{0},x_{i}]$ joining~$x_{0}$ to~$S_{i}$
in~$B$.
There exists a tree~$T$ in the union of the~$[x_{0},x_{i}]$
containing~$x_{0}$
with endpoints~$x_{i}$.

Let~$\wh X := (\overline{X} \setminus \Int B) \cup T$ and $\wh B := B
\cup (\cup_{i} C_{i})$. Notice that~$\wh X$ is (path) connected.
Indeed, every point~$x\in X\setminus \Int B$ can be connected to
some~$S_i$ by a path in~$X\setminus \Int B$ (every path from~$x_{0}$
to~$x$ intersects~$S$), while every point of each component~$S_i$ can be
connected to~$x_{0}$ by a path in $S_{i} \cup T \subset \wh X$.  By
the results of Section~\ref{prelim}, the triad $(\ov X; \wh X, \wh B)$ is
a~$CW$-triad.  Since every loop of~$\wh B$ can be deformed into a loop
of~$B$, the inclusion~$\wh B \subset \overline{X}$ induces a trivial
homomorphism of fundamental groups because of \propref{ball}.
Furthermore, the space~$\wh X \cap \wh B=T\cup(\cup C_i)$ is simply
connected.  Since~$\overline{X} = \wh X \cup \wh B$, we deduce from
Seifert--van Kampen theorem that the inclusion $\wh X \subset \overline X$ 
induces an isomorphism of fundamental groups.  Thus, the
relation~\eqref{eq:isom} leads to
\begin{equation}
\pi_{1}(\wh X) \simeq \pi_{1}(X) \simeq G.
\end{equation}
We endow each cone~$C_i$ over~$S_i$ with the round metric, \cf Appendix.
By \propref{conearea}, the area of~$C_i$ is equal to
$\frac{1}{\pi} (\length S_i)^2$.
Since the sum of the lengths of the~$S_i$'s is equal to the
length of
$S$,
the total area of~$\cup_i C_i$ is at most~$\frac{1}{\pi}
(\length
S)^2$.
The tree~$T$ is endowed with its standard metric, \ie the
length of
each of its
edges is equal to~$1$.
The metrics on~$X \setminus B$,~$\cup_i C_i$ and~$T$
induce a metric, noted~$\wh \gmetric$, on the union
$\wh X = (X \setminus B) \cup (\cup_i
C_i) \cup T$.

By construction, one has
\begin{equation}
\sys(\wh X) \geq \sys(X) = 1.
\end{equation}
Furthermore, we have
\begin{equation}
\area \wh X \leq \area X - \area B + \frac{1}{\pi}
(\length
S)^{2}.
\end{equation}
The inequality~\eqref{eq:assum} leads to
\begin{equation}
\area \wh X \leq \area X - \left( \lambda -
\frac{1}{\pi}
\right) (\length
S)^{2}.
\end{equation}
Hence,~$\sigma(\wh \gmetric) \leq \sigma(\gmetric)
\leq
\sigma(G) + \varepsilon$. Here, the first inequality holds since
$\lambda > \frac{1}{\pi}$ while the second one follows from \eqref{inequality}.

Since~$\sigma(G) \leq \area(\wh X)$ and~$\area(X)
\leq
\sigma(G) +
\varepsilon$, we also obtain
\begin{equation}
\left( \lambda - \frac{1}{\pi} \right) (\length S)^{2} <
\varepsilon.
\end{equation}
Since~$\varepsilon < 4 (\lambda - \frac{1}{\pi}) (\delta -
\varepsilon)^{2}$ and~$\delta \leq r_{0}$, we deduce that
\begin{equation}
\length S < 2 (\delta - \varepsilon) \leq 2 r_{0} - 2
\varepsilon.
\end{equation}

Now, by \lemref{lem:simple} and Proposition~\ref{lem:1compo}, every pointed 
systolic loop~\mbox{$\gamma \subset X$} at~$x_{0}$
intersects exactly two path-connected components of~$S$, say~$S_{1}$ 
and~$S_{2}$ (recall that $r_{0} \geq \delta > \frac{1}{2} \left( \sys(X,x) - 
\sys(X) \right)$, \cf \eqref{eq:delta}).
Since~$\gamma$ contains an arc of~$X \setminus 
\Int(B)$ joining~$S_{1}$ 
to~$S_{2}$, \cf Lemma~\ref{lem:simple},
we conclude that two path-connected components of~$S$
coalesce forward.
Thus, by \propref{ball} and \lemref{42},~$G$ has a positive Grushko free 
index.

\m
Now, the points~$x_{1}$ and~$x_{2}$, which are joined by a path in
$X \setminus \Int B$, are also joined to~$x_{0}$ by a unique
geodesic arc in the tree~$T$.
Identify the unique edge of the tree~$T$ which contains~$x_{1}$ with
the segment~$[0,1]$.
Set~$Y:=\wh X\setminus I \subset \wh X$, where $I=(\frac{1}{3},\frac{2}{3})$.
Since~$\wh X$ is connected and the endpoints of~$I$ are joined by
a path in~$\wh X\setminus I$, we conclude that~$Y$ is connected.
Furthermore, the space~$\wh X$, obtained by gluing back the interval~$I$ 
to~$Y$, 
is homotopy equivalent to~$Y \vee S^1$. In particular,
%
\begin{equation}
G \simeq \pi_{1}(\wh X) \simeq \pi_{1}(Y) * \Z.
\end{equation}
By uniqueness (up to isomorphism) of the Grushko decomposition~$G \simeq 
F_p*H$, we obtain
\begin{equation}
\pi_{1}(Y) \simeq F_{p-1}*H.
\end{equation}
Furthermore,~$\sigma(Y) \leq \sigma(\wh \gmetric) \leq \sigma(G) +
\varepsilon$.
In particular, we deduce that~$\sigma(F_{p-1}*H) \leq
\sigma(G) + \varepsilon$, which concludes the proof of Lemma~\ref{path}.
\end{proof}

\begin{proposition} \label{prop:dicho}
With the previous notation,
\begin{enumerate}
\item[(i)] either the Grushko free index~$p$ of~$G$ is positive, and
\begin{equation} \label{eq:p}
\sigma(G) \geq \sigma(F_{p-1}*H) - \varepsilon,
\end{equation}
\item[(ii)] or
\begin{equation}\label{ineq1}
\area B(x,r) \geq \frac{1}{4\lambda} (r-\delta)^{2}
\end{equation}
for every~$x \in X$ and every~$r \in (\delta, \frac{1}{2})$.
\end{enumerate}
\end{proposition}

\begin{proof}
By Lemma~\ref{path}, we can assume that
\begin{equation}
\area B(x,r) \leq \lambda \, (\length S(x,r))^2
\end{equation}
for every $x \in X$ and $r \in (\delta,\tfrac{1}{2})$, otherwise the claim~(i) 
holds.
Now, if~$a(r)$ and~$\ell(r)$
represent the area of~$B(x,r)$ and the length of~$S(x,r)$, respectively, then
the claim~(ii) follows from Lemma~\ref{lem:ineq} below along with the coarea 
formula.
\end{proof}

\begin{lemma} \label{lem:ineq}
Fix~$\delta \in (0,\frac{1}{2})$.
Assume that, for all~$r \in (\delta, \frac{1}{2})$,
we have
\begin{equation}
a(r) := \int_{0}^{r} \ell(s) \, ds \leq \lambda \, \ell(r)^{2}.
\end{equation}
Then, for every~$r \in (\delta, \frac{1}{2})$, we have
\begin{equation}
a(r) \geq \frac{1}{4 \lambda} (r-\delta)^{2}.
\end{equation}
\end{lemma}

\begin{proof}
The function~$\ell(r)$ is a piecewise continuous positive function by 
\propref{homeo}. So, the function~$a(r)$ is continuously differentiable for all 
but finitely many~$r$ in~$(\delta, \frac{1}{2})$.
Furthermore, $a'(r) = \ell(r)$ for all but finitely many~$r$ in~$(\delta, 
\frac{1}{2})$.  By assumption, we
have
$$
a(r) \leq \lambda \, a'(r)^2
$$
for all but finitely many~$r \in (\delta, \frac{1}{2})$. That is,
$$ \left( \sqrt{a(r)} \right)' = \frac{a'(r)}{2 \sqrt{a(r)}} \geq
\frac{1}{2\sqrt{\lambda}}.~$$
Integrating this inequality from~$\delta$ to~$r$, we get
$$
 \sqrt{a(r)} \geq \frac{1}{2\sqrt{\lambda}} (r-\delta).
$$
Hence, for every~$r \in (\delta, \frac{1}{2})$, we obtain
$$
 a(r) \geq \frac{1}{4 \lambda} (r-\delta)^{2}.
$$
\end{proof}

\section{Main results}

In this section, we extend previous results for groups of zero Grushko
free index to arbitrary finitely presentable groups.  More precisely, we
establish a uniform bound on the systolic ratio of unfree finitely
presented groups and a finiteness result for the unfree part of a
group with systolic ratio bounded away from zero.  

\begin{theorem} \label{theo:general}
Every unfree finitely presentable group~$G$ satisfies
\begin{equation} \label{eq:8/pi}
\SR(G) \leq \frac{16}{\pi}.
\end{equation}
\end{theorem}

\begin{remark}
The upper bound by~$\frac{16}{\pi}$ on the systolic ratio
in~\eqref{eq:8/pi} is not as good as the upper bound by~$4$ 
obtained in Theorem~\ref{theo:smallballs} in the zero Grushko free index  case.
\end{remark}

\begin{proof}[Proof of Theorem~$\ref{theo:general}$]
Let us prove the inequality~\eqref{eq:8/pi} by induction on the Grushko free 
index  of~$G$.
To start the induction, consider a finitely presentable group~$G$ of zero 
Grushko 
free index. Then, by \theoref{theo:smallballs},
$$
\sigma(G)\geq \frac{1}{4} > \frac{\pi}{16}.
$$

Now, assume that the inequality \eqref{eq:8/pi} holds for all finitely 
presented 
groups whose Grushko free index is less than~$p$.
Consider a finitely presentable group~$G$ with positive Grushko free index~$p$. 
The group~$G$ decomposes as $G = F_{p}*H$
where $F_{p}$ is free of rank~$p$ and~$H$ is of zero Grushko free
index. We will use the notation of Section~\ref{sec:dicho}.

If the
inequality~$\eqref{ineq1}$ holds for all~$x\in X$ and~$r\in 
(\delta,\frac{1}{2})$, then
\begin{equation}
\sigma(\gmetric) = \area X \geq \frac{1}{4\lambda} \left(
\frac{1}{2}-\delta \right)^{2}.
\end{equation}
That is,
\begin{equation}
\sigma(G) \geq \frac{1}{4\lambda} \left( \frac{1}{2}-\delta
\right)^{2} - \varepsilon.
\end{equation}
Note that the right-hand term goes to~$\frac{\pi}{16}$ as
$\delta
\rightarrow 0$,~$\lambda \rightarrow \frac{1}{\pi}$ and
$\varepsilon \rightarrow 0$. Thus,~$\sigma(G)\ge \frac{\pi}{16}$, \ie the 
inequality~$\eqref{eq:8/pi}$ holds.

Therefore, we can assume that the inequality~\eqref{eq:p} holds, \ie
\[
\sigma(G) \geq \sigma(F_{p-1}*H) - \varepsilon.
\]
By induction on~$p$, we obtain
\begin{equation} \label{peps}
\sigma(G) \geq \frac{\pi}{16} - \varepsilon.
\end{equation}
This implies the inequality \eqref{eq:8/pi} as~$\varepsilon \to 0$. 
\end{proof}

\begin{theorem}\label{theo:bound}
Let $G$ be a finitely presentable group. If~$\sigma(G) < C$ for some~$C>0$, then 
the isomorphism class of the unfree factor 
of~$G$ lies in the finite set~$\Gamma' \left( \frac{576\,C}{\pi}
\right)$.  

In particular, the number of isomorphism classes of unfree
factors of finitely presentable groups~$G$ such that~$\sigma(G)<C$ is at
most
\begin{equation} \label{eq:AC3}
A^{C^{3}},
\end{equation}
where~$A$ is an explicit universal constant.
\end{theorem}

\begin{remark}
Theorem~\ref{theo:bound} answers, to a certain extent, a
question of M.~Gromov, \cf~\cite[p.~337]{Gr2}.
\end{remark}

\begin{proof}[Proof of Theorem~$\ref{theo:bound}$]
We prove the result by induction on the Grushko free index of~$G$.
Theorem~\ref{theo:finite} shows that the isomorphism class of every finitely 
presented group~$H$ of zero Grushko 
free index with~$\sigma(H) < C$ lies in $\Gamma(144\,C) \subset 
\Gamma'(\frac{576\,C}{\pi})$.

Now, let~$G$ be a finitely presentable group of positive Grushko free index~$p$, 
that is~$G=F_p*H$
where~$H$ has zero Grushko free index.
Suppose that~$\sigma(G) < C$.
We will use the notation of Section~\ref{sec:dicho}.
Note that we can always assume that~$\sigma(\gmetric) < C$ for $\gmetric$ as in 
\eqref{inequality}.

If the inequality~$\eqref{ineq1}$ holds for all~$x\in X$ and
$r\in (\delta,\frac{1}{2})$, then the inequality~\eqref{alpha} holds for
$\alpha = \frac{1}{4\lambda} (\frac{1}{12}-\delta)^{2}$.
Hence, by \propref{prop:finite}, the isomorphism class of the group~$G$
lies in the finite set~$\Gamma(C/\alpha)$, which is contained in
$$
\Gamma\left(\frac{576\,C}\pi\right)\subset
\Gamma'\left(\frac{576\,C}\pi\right)
$$
if $\delta$ is close enough to~$0$ and $\lambda$ is close enough 
to~$\frac{1}{\pi}$.
In particular, this shows that the isomorphism class of the group~$H$ lies in
$
\Gamma'\left(\frac{576\,C}\pi\right).
$

So, we can assume that the inequality~\eqref{eq:p} holds.
Since~\mbox{$\sigma(G) < C$}, we obtain
\[
\sigma(F_{p-1}*H) < C+\varepsilon.
\]

By induction on~$p$, we derive that the isomorphism class of~$H$ lies
in~$\Gamma'\left(\frac{576\,(C+\varepsilon)}{\pi}\right)$
for all~$\varepsilon>0$.
Thus, the isomorphism class of~$H$ lies
in~$\Gamma'\left(\frac{576\,C}{\pi}\right)$

Finally, by \lemref{lem:count}, we can take
$A=2^{\left(\frac{576}{\pi}\right)^{3}}$ in~\eqref{eq:AC3}.
\end{proof}

We have the following Corollary that generalizes Example~\ref{moore}.

\begin{cory}\label{sequence}
Let~$G_1, \ldots, G_n, \ldots~$be a sequence of pairwise
non-isomorphic groups of bounded Grushko free index. Then
$$
\lim_{n \rightarrow \infty} \SR(G_{n}) = 0.
$$
\end{cory}

\begin{proof} This follows from \theoref{theo:bound} since,
given~$\varepsilon > 0$, there are only finite number of~$n$'s with
$\SR(G_{n}) <\varepsilon$.
\end{proof}

\begin{example}
Let $G_n$ be the free product of $n$ unfree finitely presentable groups.
As in Corollary~\ref{sequence}, we obtain from Theorem~\ref{theo:bound} that 
the 
systolic ratio of the sequence~$\{G_n\}$ tends to zero as $n \to \infty$, \cf 
\cite[p.~337]{Gr2}
\end{example}

\appendix
\section{Round metrics}\label{round}

Consider the upper hemisphere~$H$ of the radius~$r$,
\[
H:=\{ (x,y,z)\in \R^3 \bigm| x^2+y^2+z^2 = r^{2}, z\ge 0 \}
\]
We equip $H$ with the sphere metric $\dist_H$. Let~$K=\{(x,y,z)\in 
H\bigm|z=0\}$ and~$p=(0,0,r)\in H$. Given a point~$q\in H,
q\ne p$, consider the geodesic arc of length~$\pi r/2$ that starts at~$p$,
passes through~$q$ and ends at some point~$x=x(q)\in K$. We define~$t=t(q)$ as
the length of the geodesic segment joining~$p$ and~$q$. Clearly,~$q$ determines 
and is uniquely determined by~$x$ and~$t$. Thus, every point of~$H$ can be 
described
as a pair~$(x,t)$ where~$x\in K$ and~$t\in [0, \pi\,r/2]$. Here,~$(x,0)=p$ for 
all~$x$.

We define a function~$f: [0,\pi r]\times [0, \pi r/2]^2 \to \R$ by setting
\begin{equation}\label{f}
f(R,t_1,t_2)=\dist_H((x_1,t_1),(x_2, t_2))
\end{equation}
where~$(x_i, t_i)\in H, i=1,2$ are such that~$\dist_K(x_1,x_2)=R$. Clearly, the
function~$f$ is well-defined.

Now, let~$S$ be a finite metric graph of total length~$L$ and set~$r=L/\pi$.
Consider the cone~$C=(S\times [0,\pi\,r/2])/(S\times \{0\})$. Every point
of~$C$ can be written as~$(x,t)$ where~$x\in S$ and~$t\in [0,\pi r/2]$.
We denote by~$v$ the vertex~$(x,0)$ of the cone.
We equip~$C$ with a piecewise smooth metric by setting
\begin{equation*}
\dist_C((x_1,t_1),(x_2,t_2))=f(\dist_S(x_1,x_2),t_1,t_2)
\end{equation*}
where~$f$ is the function defined in~\eqref{f}. It is clear that~$\dist_C$ is a
metric since~$\dist_H$ is, and it is piecewise smooth since~$\dist_S$ is.
We call
this metric the {\it round metric} on~$C$. Clearly, the inclusion~$S\subset C$
is an isometry.

Furthermore, the region~$(e \times [0,\pi r/2])/(e \times \{0\})$ of~$C$, 
where~$e$ is an edge of~$S$,
is isometric to a sector of the hemisphere~$H$ of angle~$\frac{1}{r} 
\length(e)$.
Thus, the area of this region is equal to~$r \length(e)$.
We immediately deduce the following result.

\begin{prop}\label{conearea}
The area of the cone~$C$ is given by
$$
\area(C)=rL= L^2/\pi.
$$
\end{prop}

\end{document}